\documentclass[10pt,english]{article}

\usepackage{geometry}
\usepackage{verbatim,times}
\usepackage{float}
\usepackage{bm}
\usepackage{amsmath,amsthm}
\usepackage{amssymb}
\usepackage{graphicx}
\usepackage{subfig}
\usepackage{hyperref}
\hypersetup{
 colorlinks,linkcolor=red,anchorcolor=blue,citecolor=blue}
\usepackage{breakurl}
\usepackage{tikz}
\usetikzlibrary{shapes,arrows,positioning,calc}

\makeatletter

\tikzstyle{block} = [draw, thick, rectangle, 
minimum height=3em, minimum width=6em]
\tikzstyle{pt} = [coordinate]
\tikzset{
	block/.style={
		draw, 
		rectangle, 
		minimum height=0.8cm, 
		minimum width=0.8cm, 
		align=center
	}
}

\theoremstyle{plain} \newtheorem{lem}{\textbf{Lemma}} \newtheorem{prop}{\textbf{Proposition}}\newtheorem{thm}{\textbf{Theorem}}
\newtheorem{cor}{\textbf{Corollary}} 
 \newtheorem{defn}{\textbf{Definition}}
 
\theoremstyle{rem}\newtheorem{rem}{\textbf{Remark}}

\begin{document}

\title{Analytical Convergence Regions of Accelerated Gradient Descent
	in Nonconvex Optimization under Regularity Condition}
\author
{
	Huaqing Xiong\thanks{Department of Electrical and
		Computer Engineering, The Ohio State University, Columbus, OH 43210, USA; Email:
		\texttt{\{xiong.309, zhang.491\}@osu.edu}.}
	\and Yuejie Chi\thanks{Department of Electrical and Computer Engineering, Carnegie Mellon University, Pittsburgh, PA 15213, USA; Email:
		\texttt{yuejiechi@cmu.edu}. }
	\and Bin Hu\thanks{Department of Electrical and Computer Engineering and Coordinated Science Laboratory, University of Illinois at Urbana-Champaign, Urbana, IL 61801, USA; Email: \texttt{binhu7@illinois.edu}.}
	\and Wei Zhang\footnotemark[1]
}

\date{October 2018}

\maketitle

\begin{abstract}

There is a growing interest in using robust control theory to analyze and design optimization and machine learning algorithms. This paper studies a class of nonconvex optimization problems whose cost functions satisfy the so-called Regularity Condition (RC). Empirical studies show that accelerated gradient descent (AGD) algorithms (e.g. Nesterov's acceleration and Heavy-ball) with proper initializations often work well in practice. However, the convergence of such AGD algorithms is largely unknown in the literature. The main contribution of this paper is the analytical characterization of the convergence regions of AGD under RC via robust control tools. Since such optimization problems arise frequently in many applications such as phase retrieval, training of neural networks and matrix sensing, our result shows promise of robust control theory in these areas. 

\end{abstract}

\section{Introduction} \label{sec:Introduction}

Recently, control theoretic tools have gained popularity in analyzing optimization and machine learning algorithms \cite{lessard2016analysis,hu2017dissipativity,hu2017unified,pmlr-v80-hu18b,fazlyab2018analysis,wilson2016lyapunov,nishihara2015general,sundararajan2017robust,cherukuri2017role,hu2017analysis}. Typically, convergence analysis in optimization is performed in a case-by-case manner and the corresponding techniques highly depend on the structure of algorithms and assumptions of objective functions. However, by representing iterative algorithms and prior information of objective functions as feedback dynamical systems, we can apply tools from control theory to carry out the convergence analysis in a more systematic way. 

Such a framework is pioneered by \cite{lessard2016analysis}, where the authors used semidefinite programming to analyze the convergence of a class of optimization algorithms including gradient descent (GD), Heavy-ball (HB) and Nesterov's accelerated gradient (NAG), by assuming the gradient of the loss function satisfies some integral quadratic constraints (IQCs). Indeed, the standard smoothness and convexity assumptions can be well rephrased as IQCs.
Afterwards, similar approaches have been developed to analyze various optimization algorithms such as  ADMM \cite{nishihara2015general}, distributed methods \cite{sundararajan2017robust}, proximal algorithms \cite{fazlyab2018analysis}, and stochastic finite-sum methods \cite{hu2017unified,pmlr-v80-hu18b}. Exploiting the connection between control and optimization also provides new insights into momentum methods \cite{hu2017dissipativity,wilson2016lyapunov,HuACC2017} and facilitates the design of new algorithms \cite{van2018fastest,cyrus2018robust,dhingra2018proximal,pmlr-v80-kolarijani18a}. Moreover, control tools are also useful in analyzing the robustness of algorithms against computation inexactness \cite{lessard2016analysis,cherukuri2017role,hu2017analysis,aybat2018robust}.

This paper considers a class of nonconvex optimization problems whose objective functions satisfy the so-called Regularity Condition (RC)~\cite{candes2015phase,chi2019nonconvex}, which is a geometric condition characterizing the curvatures of nonconvex functions. Such a condition has appeared in many important machine learning and signal processing applications including phase retrieval \cite{candes2015phase,ChenCandes15solving,zhang2017reshaped,zhang2016provable,wang2017solving}, deep linear neural networks \cite{zhou2017characterization}, shallow nonlinear neural networks \cite{li2017convergence}, matrix sensing \cite{tu2016low,li2018median}, to name a few. While it is straightforward to show that GD converges linearly for problems satisfying RC \cite{candes2015phase}, however, the behavior of AGD under RC is remains poorly understood in theory, despite its empirical success \cite{pauwels2017fienup}.

Our work is motivated to deepen the understanding of convergence of AGD under RC. The main contribution of this paper lies in the theoretical convergence guarantee of AGD algorithms under RC. In particular, we provide an analytical
characterization of hyperparameter choices that ensure
linear convergence of AGD for all functions that satisfy RC. In addition, the framework and tools  developed herein may inspire more research on analyzing the convergence of sophisticated algorithms under nonconvex settings. Specifically, the analysis for momentum methods typically involves subtle constructions of Hamiltonians (or Lyapunov functions). Our frequency domain approach can implicitly ensure the existence of such Lyapunov functions without explicit constructions. This sheds new light on the analysis of momentum methods for nonconvex optimization problems.

It is worth noting that, RC is essentially the same as the sector bound condition in \cite{lessard2016analysis} when it holds globally. In \cite[Sections 4.5 and 4.6]{lessard2016analysis},  LMI conditions have been implemented to analyze momentum methods under the sector bound condition.  
The results, however, are numerical validations of convergence for given hyperparameters by running a small semidefinite program, and such validation needs to be done whenever the hyperparameters are changed. Built on this prior work, we focus on how to obtain analytical convergence regions from these LMIs without solving semidefinite programs. Our analytical results for momentum methods under RC provide deeper understandings on the connection between control theory and nonconvex optimization.

The rest of the paper is organized as follows. In Section 2, we state the problem and introduce the AGD methods of interest. Section 3 presents how to incorporate the algorithms and RC into a dynamical system and transfer the convergence analysis into the stability analysis. Section 4 derives the analytical convergence conditions of AGD methods under RC via a frequency domain approach by applying the KYP lemma. In Section 5, we discuss how to extend the results to the general case where RC only holds locally.

\section{Problem background}

A general optimization problem can be described as
\begin{equation}\label{eq:loss}
\underset{z\in\mathbb{R}^{n}}{\text{minimize}}\quad f(z),
\end{equation}
where $f(z)$ may be both nonconvex and nonsmooth.

\subsection{Regularity Condition}

This paper focuses on a special yet important case of nonconvex optimization, where the objective function $f(\cdot)$ satisfies the Regularity Condition (RC) \cite{candes2015phase}, defined as follows.

\begin{defn}[Regularity Condition]\label{def:RC} A function $f(\cdot)$ is said to satisfy the Regularity Condition RC($ \mu,\lambda, \epsilon $) with positive constants $ \mu,\lambda$ and $\epsilon $, if
	\begin{equation}\label{eq:locRC}
	\langle \nabla f(z), z- x^{\star} \rangle\geq\frac{\mu}{2}\| \nabla f(z)\|^2+\frac{\lambda}{2} \left\| z- x^{\star} \right\|^2
	\end{equation}
	for all $z\in \mathcal{N}_\epsilon(x^{\star}): = \left\{ z: \| z - x^{\star}\|\leq \epsilon \|x^{\star}\| \right\}$, where $x^{\star}$ is a local minimizer of $f(z)$.
\end{defn}

It is straightforward to check that one must have $\mu\lambda\leq 1$ by Cauchy-Schwartz inequality. RC can be regarded as a combination of one-point strong convexity and smoothness \cite{chi2019nonconvex}, and does not require the function $f(\cdot)$ to be convex. 
RC has appeared in a wide range of applications, a partial list of which includes phase retrieval \cite{candes2015phase}, deep linear neural networks \cite{zhou2017characterization}, shallow nonlinear neural networks \cite{li2017convergence} and matrix sensing \cite{tu2016low,li2018median}.

Our framework can handle the general case where RC only holds locally as defined. To simplify the presentation, we will first assume RC holds globally, i.e. $\epsilon=\infty$ in Definition \ref{def:RC},  in which case $x^\star$ becomes the global minimizer correspondingly. It turns out that our main results can be directly applied to the case when RC holds locally, using proper initializations, which will be explained in detail in Section \ref{sec:LRC}. Without ambiguity, we will omit the neighborhood radius $\epsilon$ and use the notation $RC(\mu,\lambda)$ to denote the global RC in the derivation of the main results.

\subsection{Accelerated gradient descent methods\label{subsec:acceleratedGD}}

In practice, AGD methods are widely adopted for its ability to accelerate the convergence. Two widely-used acceleration schemes include Nesterov's accelerated gradient (NAG) method~\cite{nesterov2003introductory}, given as
\begin{align}\label{Nesterov}
y_k &= (1+\beta)z_k-\beta z_{k-1},    \nonumber \\
z_{k+1} &= y_k-\alpha\nabla f(y_k), \quad k=0,1,\ldots,  
\end{align}
where $\alpha>0$ is the step size, $0\leq\beta<1$ is the momentum parameter; and Heavy-Ball (HB) method \cite{polyak1964some}, given as 
\begin{align}\label{HB}
y_k &= (1+\beta)z_k-\beta z_{k-1},    \nonumber \\
z_{k+1} &= y_k-\alpha\nabla f(z_k), \quad k=0,1,\ldots,
\end{align} 
where $\alpha>0$ is the step size, $0\leq\beta<1$ is the momentum parameter. In fact, we can describe a general AGD method that subsumes HB and NAG as special cases:
\begin{equation}\label{accAlg}
\begin{aligned}
y_k &= (1+\beta_2)z_{k}-\beta_2 z_{k-1},\\
z_{k+1} &= (1+\beta_1)z_{k}-\beta_1 z_{k-1}-\alpha \nabla  f(y_k).
\end{aligned}
\end{equation}

Despite the empirical success, the convergence of AGD in the nonconvex setting remains unclear to a large extent. For example, it is not known whether AGD converges under RC, whether it converges linearly if it does and how to set the step size and the momentum parameters to guarantee its convergence. These challenging questions motivate us to look for new tools to better understand AGD for nonconvex optimization.

\section{A control view on the convergence analysis of AGD under RC}

Robust control theory has been tailored to the convergence analysis of optimization methods \cite{lessard2016analysis,nishihara2015general,hu2017dissipativity,hu2017unified}.
The proofs of our main theorems also rely on such techniques. 
In this section, we will discuss how to transform the convergence analysis of AGD under RC to the robust stability analysis of dynamical systems and derive LMI conditions to guarantee the convergence.

\subsection{AGD as feedback dynamical system}
Observe that a general AGD method \eqref{accAlg} can be viewed as a linear dynamical system subject to nonlinear feedback:
\begin{equation}\label{dysys}
\begin{aligned}
z_{k+1}^{(1)} &= (1+\beta_1)z_{k}^{(1)}-\beta_1 z_{k}^{(2)}-\alpha u_k,\\
z_{k+1}^{(2)} &= z_{k}^{(1)},\\
y_k &= (1+\beta_2)z_{k}^{(1)}-\beta_2 z_{k}^{(2)},\\
u_k &= \nabla f(y_k).
\end{aligned}
\end{equation}
To see this, let $z_{k}^{(1)}=z_{k}$, and  $z_{k}^{(2)}=z_{k-1}$. Then it can be easily verified that \eqref{dysys} represents HB when $(\beta_1,\beta_2)=(\beta,0)$, and NAG when $(\beta_1,\beta_2)=(\beta,\beta)$.

\begin{figure}
	\centering
	\tikzstyle{thickarrow}=[line width=2mm,draw=osu,-triangle 45]
	\begin{tikzpicture}[scale=2]
	\node[block] (sys) {G};
	\node[block, below=1em of sys] (feedback) {$ \nabla f(y_k) $};	
	\draw[->,thick] (sys.east) -- ++ (2em,0) coordinate[yshift=1em](l){} |-node[near start, right]{\small $y_k$} (feedback.east); 
	\draw[->,thick] (feedback.west) -- ++ (-1.5em,0) coordinate[yshift=-1em](l){} |-node[near start, left]{\small $u_k$} (sys.west); 
	\end{tikzpicture}
	\caption{\small The dynamical system representation of first-order methods.}
	\label{sysdiag}
\end{figure}
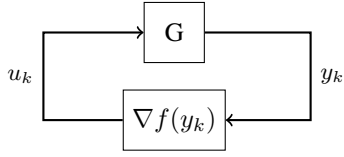

Let $\otimes$ denote the Kronecker product.
We use the notation $G(A, B, C, D)$ to denote a dynamical system $G$ governed by the following iterative state-space model:
\begin{align*}
\phi_{k+1}&=(A\otimes I_n) \phi_k+(B\otimes I_n) u_k,\\
y_k&=(C\otimes I_n) \phi_k+(D\otimes I_n) u_k,
\end{align*}
where $I_n$ is the identity matrix of size $n$. If we define $\phi_k=\left[\begin{array}{c} z_{k}^{(1)}\\z_{k}^{(2)} \end{array}\right]$ as the state, $u_k$ as the input and $y_k$ as the output, then \eqref{dysys} can be regarded as a dynamical system shown in Figure~\ref{sysdiag}, where the feedback $\nabla f(y_k)$ is a static nonlinearity that depends on the gradient of the loss function, and $G(A, B, C, D)$ is a linear system specified by 
\begin{equation}\label{eq:unisys}
\left[\begin{array}{c|c}
A  &B\\
\hline
C  &D \\
\end{array}
\right] =
\left[\begin{array}{cc|c}
1+\beta_1 &-\beta_1  &-\alpha \\
1 & 0 & 0\\
\hline
1+\beta_2 &-\beta_2  & 0
\end{array}
\right].
\end{equation}

\subsection{Convergence analysis of AGD under RC via stability analysis of a feedback system}\label{subsec:lyFunc}

In the following, we will illustrate how to show the convergence of AGD. First, define $\phi_*=\left[\begin{array}{c} x^\star\\x^\star \end{array}\right]$ as the equilibrium of the dynamical system \eqref{dysys}. If the system is (globally) asymptotically stable, then $\phi_k \xrightarrow{k\rightarrow\infty}\phi_*$. It further implies $z_{k} \xrightarrow{k\rightarrow\infty}x^\star$. In other words, the asymptotic stability of the dynamical system can indicate the convergence of the iterates to a fixed point. From now on, we  focus on the stability analysis of the feedback control system \eqref{dysys}.

The stability analysis of the feedback control system \eqref{dysys} can be carried out using robust control theory.
The main challenge is on the nonlinear feedback term $u_k = \nabla f(y_k)$. Our key observation is that RC can be rewritten as the following quadratic constraint:
\begin{equation}\label{quadrabound}
\left[ \begin{array}{c}
y_k-y_*\\u_k-u_*\end{array} 
\right]^T \left[ \begin{array}{c c}
-\lambda I_n  &I_n\\
I_n  &-\mu I_n
\end{array} 
\right] 
\left[ \begin{array}{c}
y_k-y_*\\u_k-u_*\end{array} 
\right] \geq  0,
\end{equation}
where $y_*=x^\star$ and $u_*=\nabla f(y_*)=0$.
Applying the quadratic constraint  framework in \cite{lessard2016analysis}, we can derive an LMI as a sufficient stability condition as stated in the following proposition. A formal proof of this result can be found in the appendix. 

\begin{prop}\label{thm:convergence}
	Let $x^\star\in\mathbb{R}^n$ be the global minimizer of the loss function $f(\cdot)$ which satisfies RC($ \mu,\lambda $). For a given first-order method characterized by $G(A,B,C,D)$, if there exists a matrix $P\succ0$ and $\rho\in (0,1)$ such that the following inear matrix inequality (LMI) \eqref{eq:maintheo} holds,
	\begin{equation}\label{eq:maintheo}
	\left[\begin{array}{cc}
	A^TPA-\rho^2 P &A^TPB \\
	B^TPA &B^TPB
	\end{array}
	\right]+  
	\left[ \begin{array}{c c}
	C & D\\
	0_{1\times 2} &1
	\end{array} 
	\right]^T\left[ \begin{array}{c c}
	-\lambda &1\\
	1 &-\mu
	\end{array} 
	\right]
	\left[ \begin{array}{c c}
	C & D\\
	0_{1\times 2} &1
	\end{array} 
	\right]\preceq0,
	\end{equation}
	then the state $\phi_k$ generated by the first-order algorithm $G(A,B,C,D)$ converges to the fixed point $\phi_*$ linearly, i.e.,
	\begin{equation} \label{eq:linear_convergence}
	\|\phi_k-\phi_*\|\leq \sqrt{\mathrm{cond}(P)}\rho^k\|\phi_0-\phi_*\| \text{ for all } k, 
	\end{equation}
	where $\mathrm{cond}(P)$ is the condition number of $P$.
\end{prop}

\begin{rem}
	For fixed $(A,B,C,D)$ and $\rho$, the LMI \eqref{eq:maintheo} is linear in $P$ and hence an LMI. The size of this LMI is $3 \times 3$, and the decision variable $P$ is a $2\times 2$ matrix. The size of the LMI \eqref{eq:maintheo} is independent of the state dimension $n$. 
\end{rem}

\begin{rem}
	The LMI \eqref{eq:maintheo} is similar to the one derived in~\cite{lessard2016analysis} under the so-called sector bound condition. The relationship between RC and the sector bound is discussed in detail in the appendix. Different from~\cite{lessard2016analysis}, we focus on deriving analytical convergence regions (see Section \ref{sec:anaResultMain}), in contrast to verifying convergence numerically for specific parameters, offering deeper insight regarding the convergence behavior of AGD methods under RC. In addition, we also extend the results to the case where RC holds only locally around the fixed point (see Section \ref{sec:LRC}). 
\end{rem}

\section{Convergence conditions of AGD}\label{sec:anaResultMain}
In this section, we focus on how to obtain analytical convergence conditions of AGD under RC based on \eqref{eq:maintheo}. 

Analytically solving \eqref{eq:maintheo} is challenging since one typically needs to express $P$ explicitly as a function of $(A,B,C,D)$ and $(\lambda, \mu)$. Our main idea is to transform the LMI \eqref{eq:maintheo} to equivalent frequency domain inequalities (FDIs) which can reduce unknown parameters using the classical KYP lemma~\cite{rantzer1996kalman}. Then we can derive the main convergence results by solving the FDIs analytically.

\subsection{The Kalman-Yakubovich-Popov (KYP) lemma}

We first introduce the KYP lemma and the reader is referred to \cite{rantzer1996kalman} for an elegant proof.

\begin{lem}
	(\cite[Theorem 2]{rantzer1996kalman}) Given $A$, $B$, $M$, with $\text{det}(e^{j\omega}I-A)\neq 0 $ for $\omega\in\mathbb{R}$, the following two statements are equivalent:
	\begin{enumerate}
		\vspace*{-0.5em}
		\item $\forall \omega\in\mathbb{R}$,
			\begin{equation}
			\hspace*{-0.2cm}\left[\!\begin{array}{c}
			(e^{j\omega}I-A)^{-1}B\\I
			\end{array}\!
			\right]^*\!M\!\left[\!\begin{array}{c}
			(e^{j\omega}I-A)^{-1}B\\I
			\end{array}\!
			\right]\!\prec\! 0.
			\end{equation}
		\item There exists a matrix $P\!\in\!\mathbb{R}^{n\times n}$ such that $P\!=\!P^T$ and
			\begin{equation}\label{LMIK}
			M+\left[\begin{array}{cc}
			A^TPA-P &A^TPB \\
			B^TPA &B^TPB
			\end{array}
			\right]\prec 0.
			\end{equation} 
	\end{enumerate}
\end{lem}

The general KYP lemma only asks $P$ to be symmetric instead of being positive definite (PD) as in our problem. To ensure that the KYP lemma can be applied to solve \eqref{eq:maintheo}, some adjustments of the lemma are necessary. In fact, we observe that if $A$ of the dynamical system is Schur stable and the upper left corner of $M$, denoted as $M_{11}$, is positive semidefinite (PSD), then by checking the principal minor $M_{11}+A^TPA-P\prec0$, we know $P$ satisfying \eqref{LMIK} must be PD. We define these conditions on $A$ and $M$ as KYP Conditions, which are restrictions to make sure that all solutions of symmetric $P$ for \eqref{LMIK} are PD.
\begin{defn}[KYP Conditions] The KYPC($A,M$) are listed as:
	\vspace*{-0.5em}
	\begin{enumerate}
		\item $\text{det}(e^{j\omega}I-A)\neq 0 $ for $\omega\in\mathbb{R}$;
		\item $A$ is Schur stable;
		\item The left upper corner of $M$ in \eqref{LMIK} is PSD.
	\end{enumerate}	
\end{defn}

Thus we can conclude the following corollary.
\begin{cor}\label{cor:kyp}
	Under KYPC($A,M$), the following two statements are equivalent:
	\begin{enumerate}
		\vspace*{-0.5em}
		\item $\forall \omega\in\mathbb{R}$,
			\begin{equation}
			\hspace*{-0.2cm}\left[\!\begin{array}{c}
			(e^{j\omega}I-A)^{-1}B\\I
			\end{array}\!
			\right]^*\!M\!\left[\!\begin{array}{c}
			(e^{j\omega}I-A)^{-1}B\\I
			\end{array}\!
			\right]\!\prec\! 0.
			\end{equation}
		\item There exists a matrix $P\in\mathbb{R}^{n\times n}$ such that $P\succ0$ and
			\begin{equation}\label{LMIKYP}
			M+\left[\begin{array}{cc}
			A^TPA-P &A^TPB \\
			B^TPA &B^TPB
			\end{array}
			\right]\prec 0.
			\end{equation} 
	\end{enumerate}
\end{cor}
One can easily check, however, that $A$ and $M$ of a general AGD in \eqref{eq:maintheo} do not satisfy the KYPC($A,M$). Therefore, we need to rewrite the dynamical system \eqref{dysys} in a different way to satisfy the KYPC($A,M$), so that its stability analysis can be done by combining Proposition \ref{thm:convergence} and Corollary \ref{cor:kyp}. In the following, we first introduce a way to rewrite the dynamical system to satisfy the KYPC($A,M$).

\subsection{How to satisfy the KYP Conditions?}

Recall that a general AGD can be written as \eqref{accAlg}.
Here we introduce a slack variable $\delta$ to rewrite the algorithm:
\begin{equation}\label{accAlgDelta}
\begin{aligned}
z_{k+1}\! = &\left(1+\delta+\beta_1+\delta\beta_2\right) z_{k}\!-\!(\beta_1+\delta\beta_2) z_{k-1}\\
&\quad -\alpha \nabla  f(y_k)-\delta y_k,\\
y_k\! = &(1+\beta_2)z_{k}-\beta_2 z_{k-1}.
\end{aligned}
\end{equation}
Observe that for any value of $\delta$, \eqref{accAlgDelta} provides the same update rule as \eqref{accAlg}. It can be viewed as a generalized representation of the dynamical systems corresponding to the targeted AGD. Similar to \eqref{dysys}, we rewrite \eqref{accAlgDelta} as a dynamical system $G(A',B',C',D')$:
	\begin{equation}\label{systemshift}
	\hspace*{-0.2cm}\begin{aligned}
	z_{k+1}^{(1)}\! &=\! \left(1\!+\!\beta_1\!+\!\delta\!+\!\delta\beta_2\right)z_{k}^{(1)}\! -\! (\beta_1\!+\!\delta\beta_2) z_{k}^{(2)}\! +\! u_k,\\
	z_{k+1}^{(2)}\! &= \!z_{k}^{(1)},\\
	y_k\! &= \!(1+\beta_2)z_{k}^{(1)}-\beta_2 z_{k}^{(2)},\\
	u_k\! &= \!-\alpha\nabla  f(y_k) - \delta y_k.
	\end{aligned}
	\end{equation}
Correspondingly,
	$$\left[\begin{array}{c|c}
	A' &B'\\
	\hline
	C' &D'\\
	\end{array}
	\right] =
	\left[\begin{array}{cc|c}
	1+\beta_1+\delta+\delta\beta_2 &-(\beta_1+\delta\beta_2) &1\\
	1 &0 &0\\
	\hline
	1+\beta_2 &-\beta_2 &0
	\end{array}
	\right].
	$$
In addition to the adjustment of the dynamics, the feedback of $G(A',B',C',D')$ also differs from that in \eqref{dysys}. As a consequence, the quadratic bound for the new feedback $u_k = -\alpha\nabla  f(y_k) - \delta y_k$ is shifted as stated in the following lemma.

\begin{lem}\label{lem:feedbshift}
	Let $f$ be a loss function which satisfies RC($ \mu,\lambda $) and $y_*=x^{\star}$ be a minimizer.
	If $u_k = -\alpha\nabla  f(y_k) - \delta y_k$, then $y_k$ and $u_k$ can be quadratically bounded as
	\begin{equation}\label{RCshift}
	\left[\begin{array}{c}
	y_{k}-y_*\\u_k-u_*\end{array}
	\right]^T M' \left[\begin{array}{c}
	y_{k}-y_*\\u_k-u_*\end{array}
	\right]\geq0.
	\end{equation}
	where $M'=\left[\begin{array}{cc}
	-\left(2\alpha\delta+\lambda\alpha^2+\mu\delta^2\right) &-\alpha-\mu\delta\\
	-\alpha-\mu\delta &-\mu
	\end{array}\right]$.
\end{lem}

Now we have general representations of $A'$, $M'$ with one unknown parameter $\delta$. We need to certify the region of $(\alpha,\beta_1,\beta_2)$ such that its corresponding ($A',M'$) has at least one $\delta$ satisfying the KYPC($A',M'$). 
\begin{lem}\label{lem:kypcalg}
	Let $f$ be a loss function which satisfies RC($ \mu,\lambda $). There is at least one representation of the dynamical system \eqref{systemshift} satisfying KYPC($A',M'$), if and only if the step size $\alpha$ and the momentum parameters $\beta_1,\beta_2$ obey the following restriction: 
	\begin{equation}\label{HbshiftLem}
	0<\alpha<\frac{2(1+\beta_1)(1+\sqrt{1-\mu\lambda})}{\lambda(1+2\beta_2)}.
	\end{equation}
\end{lem}
By Lemma \ref{lem:kypcalg}, if the parameters of a fixed AGD with $(\alpha,\beta_1,\beta_2)$ satisfy \eqref{HbshiftLem}, then all feasible symmetric $P$'s for \eqref{LMIKYP} can be guaranteed to be PD. Now we are ready to use the KYP lemma to complete the convergence analysis of an accelerated algorithm.

\subsection{Stability region of AGD under RC}

By Proposition \ref{thm:convergence}, we can solve the stability of the new system \eqref{systemshift} by finding some feasible $P\succ 0$ to the key LMI~\eqref{eq:maintheo} with respect to the corresponding $(A',B',C',D',M')$ for some rate $0<\rho<1$.

We are interested in obtaining the analytical region of $(\alpha, \beta_1, \beta_2)$ that guarantees the linear convergence of AGD under RC. We use the following strict matrix inequality without caring about a specific rate $\rho$,
\begin{equation}\label{LMIshift}
 \left[\begin{array}{cc}
A'^TPA'- P &A'^TPB' \\
B'^TPA' &B'^TPB'
\end{array}
\right]+ \left[\begin{array}{c c}
C' &0\\
0_{1\times 2} &1
\end{array}
\right]^TM'
\left[\begin{array}{c c}
C' &0\\
0_{1\times 2} &1
\end{array}
\right]\prec0. 
\end{equation}

\begin{figure*}[ht]
	\centering
	\captionsetup[subfigure]{labelformat=empty}
	\subfloat[\small (a) Fixing $\lambda=0.5$ and varying $\mu$]{{\includegraphics[width=6cm]{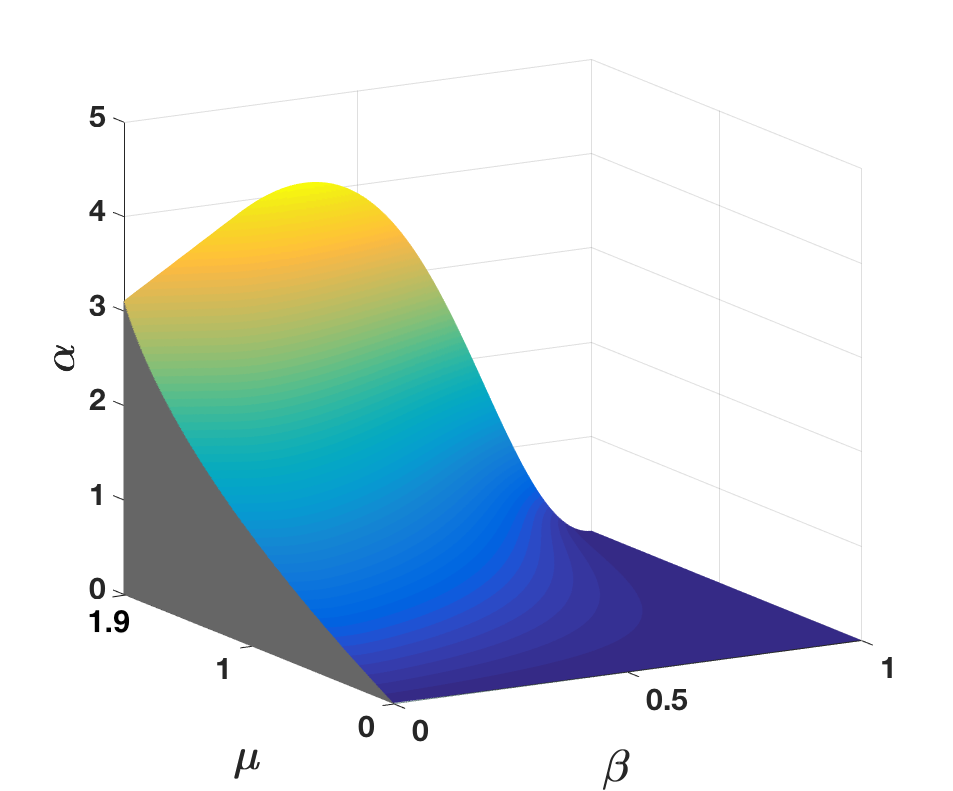} }}%
	\qquad
	\subfloat[\small (b) Fixing $\mu=0.5$ and varying $\lambda$]{{\includegraphics[width=6cm]{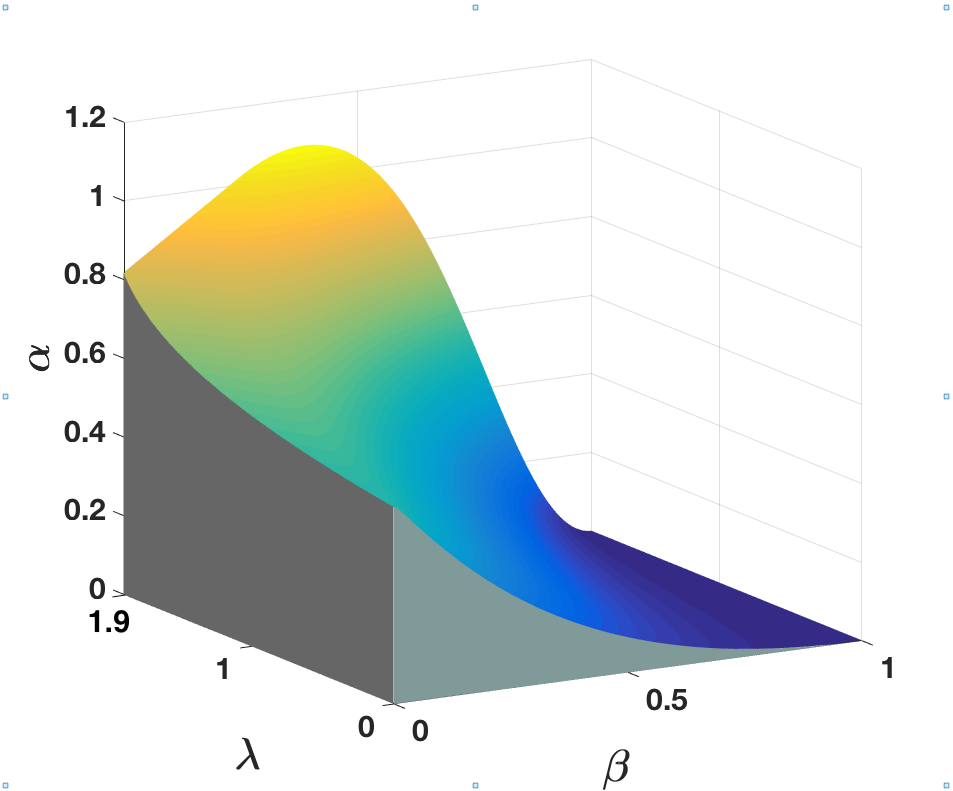} }}%
	\caption{\small Visualization of the convergence regions of HB when perturbing the RC parameters.}%
	\label{regionplot_3D}
\end{figure*}

	\begin{rem}
		Our arguments can also be modified to derive the parameter region which guarantees the convergence with a fixed rate $\rho$. For such an analysis, we can modify the LMI \eqref{LMIshift} by rescaling the matrices$(A',B')$ as $\tilde{A}=A'/\rho$ and $\tilde{B}=B'/\rho$. Then the resultant LMI can be converted to an FDI by the KYP lemma, and a similar analysis can be carried forward.  Such analytical analysis of the convergence rate is even more difficult to interpret due to the presence of $\rho$ in $(\tilde{A},\tilde{B})$. For simplicity, this paper focuses on the derivation of stability regions.
	\end{rem}

Observe that now \eqref{LMIshift} is of the same form as \eqref{LMIKYP}. By the KYP lemma (Corollary \ref{cor:kyp}), under KYPC($A',M'$) \eqref{LMIshift} can be equivalently solved by studying the following FDI:
\begin{equation}\label{FDIshift}
\left[\begin{array}{c}
(e^{j\omega}I-A')^{-1}B'\\I
\end{array}
\right]^*\left[\begin{array}{c c}
C' &0\\
0_{1\times 2} &1
\end{array}
\right]^TM' \left[\begin{array}{c c}
C' &0\\
0_{1\times 2} &1
\end{array}
\right]\left[\begin{array}{c}
(e^{j\omega}I-A')^{-1}B'\\I
\end{array}
\right]< 0,\quad \forall \omega\in\mathbb{R}.
\end{equation}

By simplifying \eqref{FDIshift} we observe that all uncertain terms can be canceled out and then conclude the following lemma.

\begin{lem}\label{lem:fdiSimp}
	To find the stability region of a general AGD method under RC($ \mu,\lambda $), or equivalently, to find the region of $(\alpha,\beta_1,\beta_2)$ such that there exists a feasible $P\!\succ\! 0$ satisfying \eqref{LMIshift}, it is equivalent to find $(\alpha,\beta_1,\beta_2)$ which simultaneously obeys \eqref{HbshiftLem} and guarantees the following FDI:
	\begin{equation}\label{FDIgeneral}
	\hspace*{-0.35cm}\begin{aligned}
	&\!4(\alpha\beta_2\!-\!\mu\beta_1)\cos^2\omega\! +\! 
	2\!\left[\! \mu(1\!+\!\beta_1)^2\!+\!\lambda\alpha^2\beta_2(1\!+\!\beta_2)\right.\!\\
	&\!\left.-\alpha(1\!+\!\beta_1)(1\!+\!2\beta_2) \right]\cos\omega\! +\!2\alpha( 1\!+\!\beta_1\!+\!2\beta_1\beta_2 )\!\\
	&\!-2\mu(1\!+\!\beta_1^2)\!-\!\lambda\alpha^2\left[ \beta_2^2\!+\!(1\!+\!\beta_2)^2 \right]< 0, \quad \forall \omega\in\mathbb{R}.\!
	\end{aligned}\!
	\end{equation}
\end{lem}
We omit the proof of Lemma \ref{lem:fdiSimp} since it follows easily from Corollary \ref{cor:kyp} and some simple calculations to simplify \eqref{FDIshift}. 

 By setting different $\beta_1,\beta_2$ in \eqref{FDIgeneral}, we can obtain the convergence condition of a general AGD method using the KYP lemma. In the following, we focus on the two most important cases: HB and NAG and other cases can be discussed in a similar way. 
The stability regions of HB and NAG can be obtained by letting $\beta_2=0$ and $\beta_1=\beta_2=\beta$ in \eqref{FDIgeneral}, respectively. Then we can obtain Theorem \ref{thm:HBregionInf} and Theorem \ref{thm:NesregionInf}.

\begin{thm}\label{thm:HBregionInf}
	Let $x^\star\in\mathbb{R}^n$ be the global minimizer of the loss function $f(\cdot)$ which satisfies RC($ \mu,\lambda $). For any step size $\alpha>0$ and momentum parameter $\beta\in(0,1)$ lying in the region: 
	$$ 
	\Big\{(\alpha,\beta): H_1(\beta)\leq\alpha \nonumber \leq\frac{2(\beta+1)(1-\sqrt{1-\mu\lambda})}{\lambda}\Big\}  \cup \Big\{(\alpha,\beta): 0<\alpha\leq\min \{H_1(\beta), H_2(\beta)
	\} \Big\}. $$
	where $H_1(\beta)  = \frac{\mu \beta^2\!+\!6\mu\beta\!+\!\mu}{\beta+1}$ and
	\begin{align*}
	H_2(\beta) & =\frac{P_2(\beta)\!-\!\sqrt{P_2(\beta)^2\!-\!4P_1(\beta)P_3(\beta)}}{2P_1(\beta)}
	\end{align*}
	with $P_1(\beta)=4\mu\lambda\beta-\beta^2-1-2\beta$, $P_2(\beta)= 2\mu\beta+2\mu\beta^2-2\mu\beta^3-2\mu$, and $P_3(\beta)=4\mu^2\beta^3+4\mu^2\beta-6\mu^2\beta^2-\mu^2\beta^4-\mu^2$,
	the iterates $z_{k}$ generated by HB~\eqref{HB} converge linearly to $x^\star$ as $k\rightarrow\infty$.
\end{thm}

\begin{thm}\label{thm:NesregionInf}
	Let $x^\star\in\mathbb{R}^n$ be {blue}the global minimizer of the loss function $f(\cdot)$ which satisfies RC($ \mu,\lambda $). For any step size $\alpha>0$ and momentum parameter $\beta\in(0,1)$ lying in the region:
	$$  
	\left\{(\alpha,\beta):N_1(\beta)\leq\alpha<\frac{2(\beta+1)(1-\sqrt{1-\mu\lambda})}{\lambda(1+2\beta)} \right\} 
	\cup \Big\{(\alpha,\beta):0<\alpha\leq\min \left\{N_1(\beta), N_2(\beta)
	\right\} \Big\}.
	$$
	where 
	$$
	N_1(\beta)=\frac{Q_1(\beta)-\sqrt{Q_1(\beta)^2-(1+6\beta+\beta^2)Q_2(\beta)}}{2\lambda\beta(\beta+1)},$$
	$$ 
	N_2(\beta)=\left\{\! \beta: \frac{Q_3(\beta)\!-\!\sqrt{Q_3(\beta)^2\!-\!(1-\beta)^2Q_2(\beta)}}{2\lambda\beta(\beta+1)}\!\leq\alpha, \; g\left(\frac{(\mu-\alpha)(1+\beta)^2+(\lambda\alpha^2-\alpha)(\beta+\beta^2)}{4\mu\beta-4\alpha\beta}\right)\!=\!0\!\right\},
	$$
	$Q_1(\beta)=1+7\beta+2\beta^2$,$Q_2(\beta)=4\mu\lambda\beta(1+\beta)$,$Q_3(\beta)=1-\beta+2\beta^2$,
	the iterates $z_{k}$ generated by NAG~\eqref{Nesterov} converge linearly to $x^\star$ as $k\rightarrow\infty$.
\end{thm}
\begin{rem}\label{rem:Nes}
	The bound $N_2(\beta)$ is an implicit function of $\beta$. It is hard to derive an explicit expression since $g(\cdot)=0$ is a 4th-order equation of $\beta$. The function $g(\eta)$ is:
	$$ 
	g(\eta):=4\mu\beta \eta^2-2(2\mu\beta+\mu\beta^2-\alpha\beta+\mu-\alpha)\eta
	+2\mu+2\mu\beta^2-2\alpha-2\alpha\beta+\lambda\alpha^2.
	$$
\end{rem}
\begin{figure*}[ht]
	\centering
	\captionsetup[subfigure]{labelformat=empty}
	\subfloat[\small (a) HB]{{\includegraphics[width=6cm]{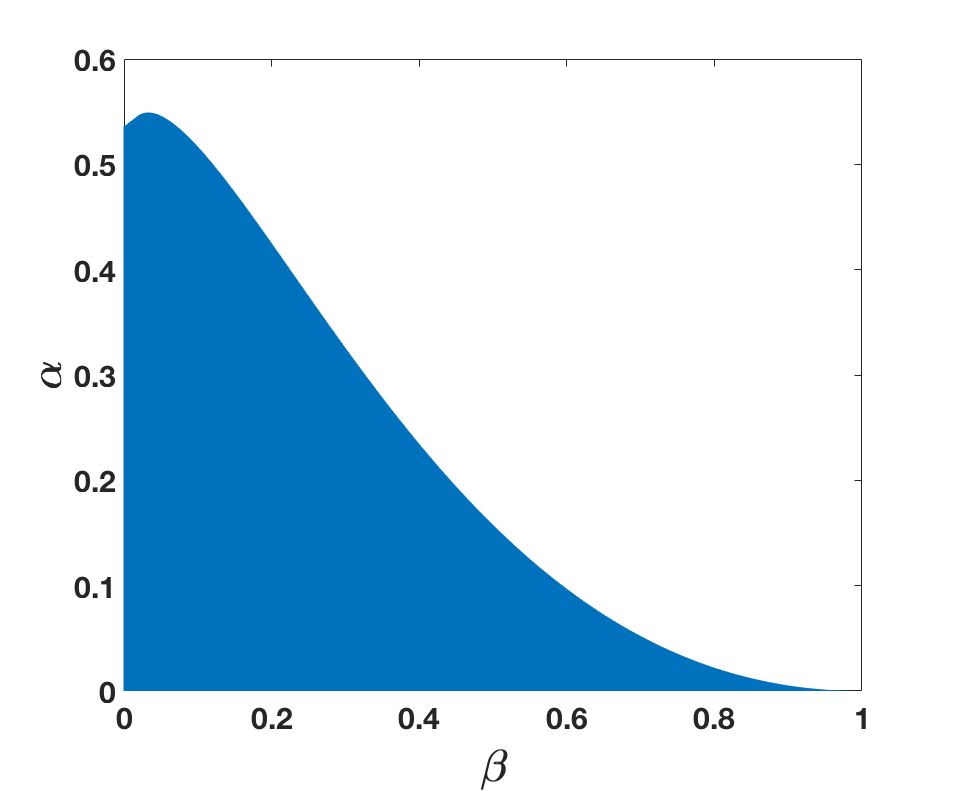} }}%
	\qquad
	\subfloat[\small (b) NAG]{{\includegraphics[width=6cm]{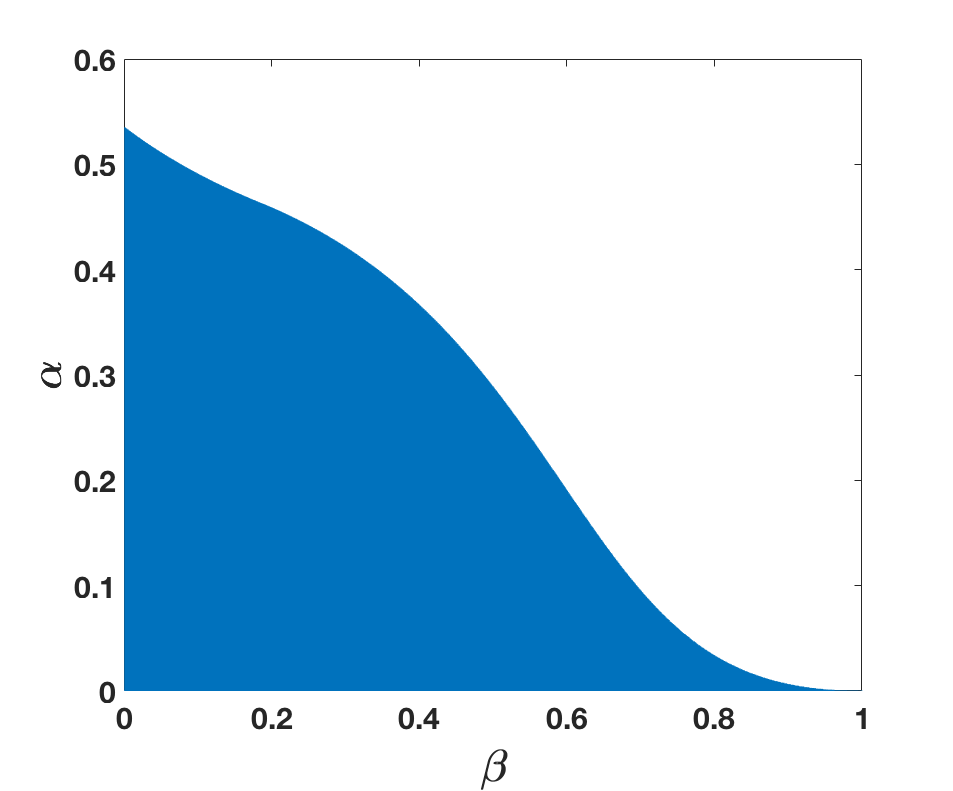} }}%
	\caption{\small Visualization of the convergence regions of two AGD methods taking RC parameters as $\mu=0.5,\lambda=0.5$.}%
	\label{regionplot}
\end{figure*}
The analytical results stated in the above theorems can provide rich insights on the convergence behaviors of AGD. Take the convergence region of HB as an example.
In Figure \ref{regionplot_3D} (a), we fix the RC parameter $\lambda=0.5$ and vary $\mu$ within $[0.01,1.9]$, while in Figure \ref{regionplot_3D} (b), we fix $\mu=0.5$ and vary the value of $\lambda$ within $[0.01,1.9]$. Observe that when we fix one of the RC parameter and increase the other, the stability region of ($\alpha,\beta$) gets larger.

Notice that $\mu$ plays a role similar to the inverse of the smoothness parameter, and therefore, it dominantly determines the step size, which is clearly demonstrated in Figure~\ref{regionplot_3D} (a). In addition, when we fix the values of a pair of ($\mu,\lambda$) (e.g. Figure~\ref{regionplot}), we can see that even when $\alpha$ exceeds the value of the bound of GD (the maximal feasible $\alpha$ when $\beta=0$), the Heavy-ball method can still ensure convergence when we choose $\beta$ properly. This property has not been discussed in the literature.

We emphasize that our theoretic analysis is a complement rather than a replacement for the numerical LMI approach in \cite{lessard2016analysis}.
Our closed-form expressions for the stability region do provide some complementary benefits to the numerical approach in \cite{lessard2016analysis}.
First, from our closed-form formulas, one can tell that the stability region of HB is well described by the feasible set of some relatively simple quadratic inequalities while the characterization of the stability region boundary of NAG partially involves a fourth-order polynomial. Such a difference is not directly reflected by the numerical LMI approach in \cite{lessard2016analysis}.
Actually our closed-form expression for the stability region of HB is quite simple. Our study on HB and NAG just illustrates that the interpretability of the analytical formulas for the stability region depends on the specific momentum method being analyzed.
Second, the stability region is easier and faster to visualize from analytical forms than numerical results. When given a pair of $(\mu,\lambda)$, one needs to make a small grid of $(\alpha,\beta_1,\beta_2)$ and solve an LMI for each single pair, which is computationally complex but can be avoided if we have closed-form analytical results. 
More importantly, the LMI conditions in Lessard et al. [16] can only certify convergence numerically for fixed $(\alpha,\beta_1,\beta_2)$ under a given pair of RC parameter $(\mu,\lambda)$. However, our analytical results provide continuous stability regions with respect to $(\mu,\lambda)$, which is hard to achieve using numerical results.

\subsection{Numerical example}

In this subsection, we will use a simple example satisfying RC globally to show show how our results help to choose parameters of different first-order methods in practice.

\begin{figure*}[ht]
	\centering
	\captionsetup[subfigure]{labelformat=empty}
	\subfloat[\small (a) Loss function $f$]{{\includegraphics[width=6cm]{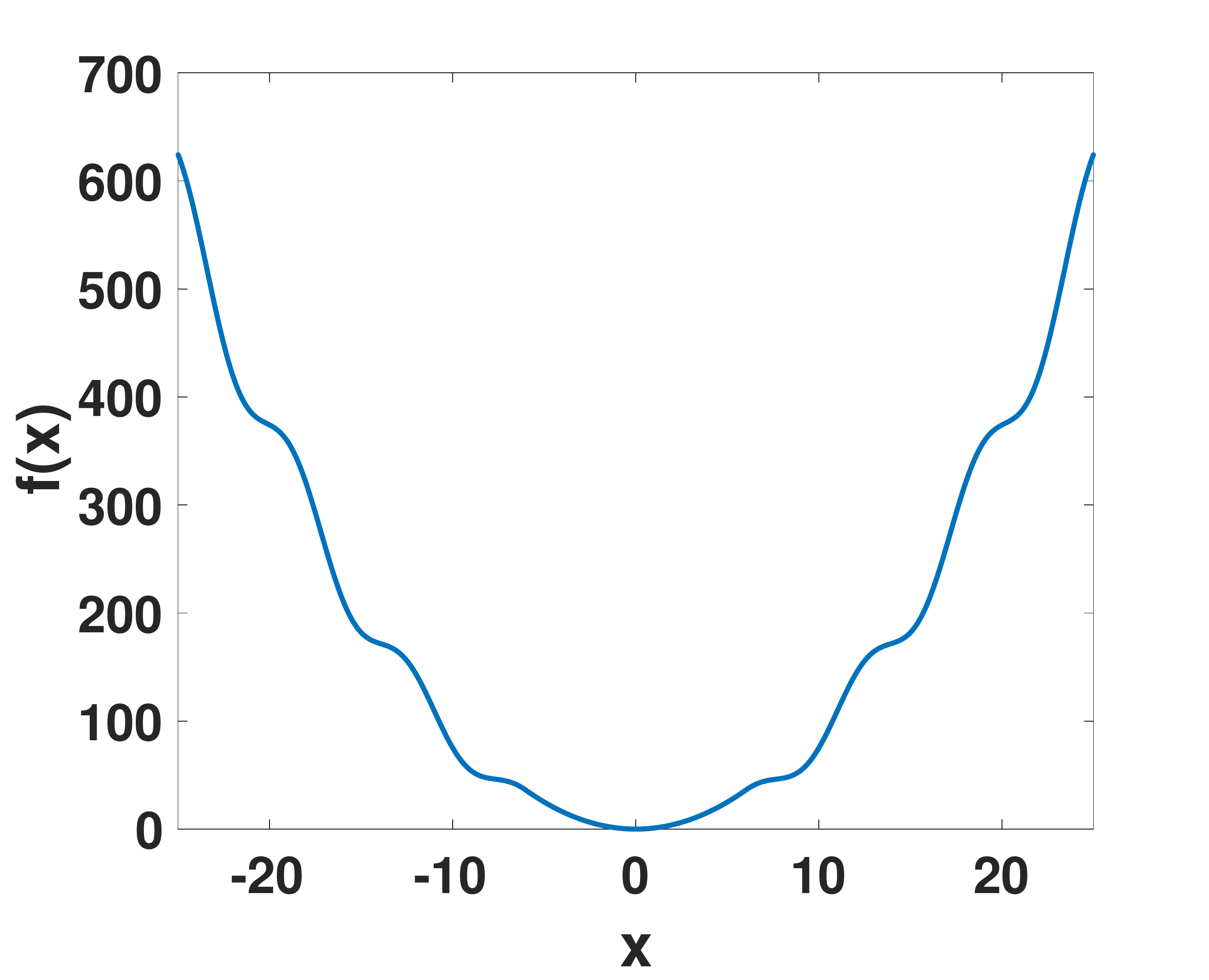} }}%
	\qquad
	\subfloat[\small (b) Convergence of three algorithms]{{\includegraphics[width=6.18cm]{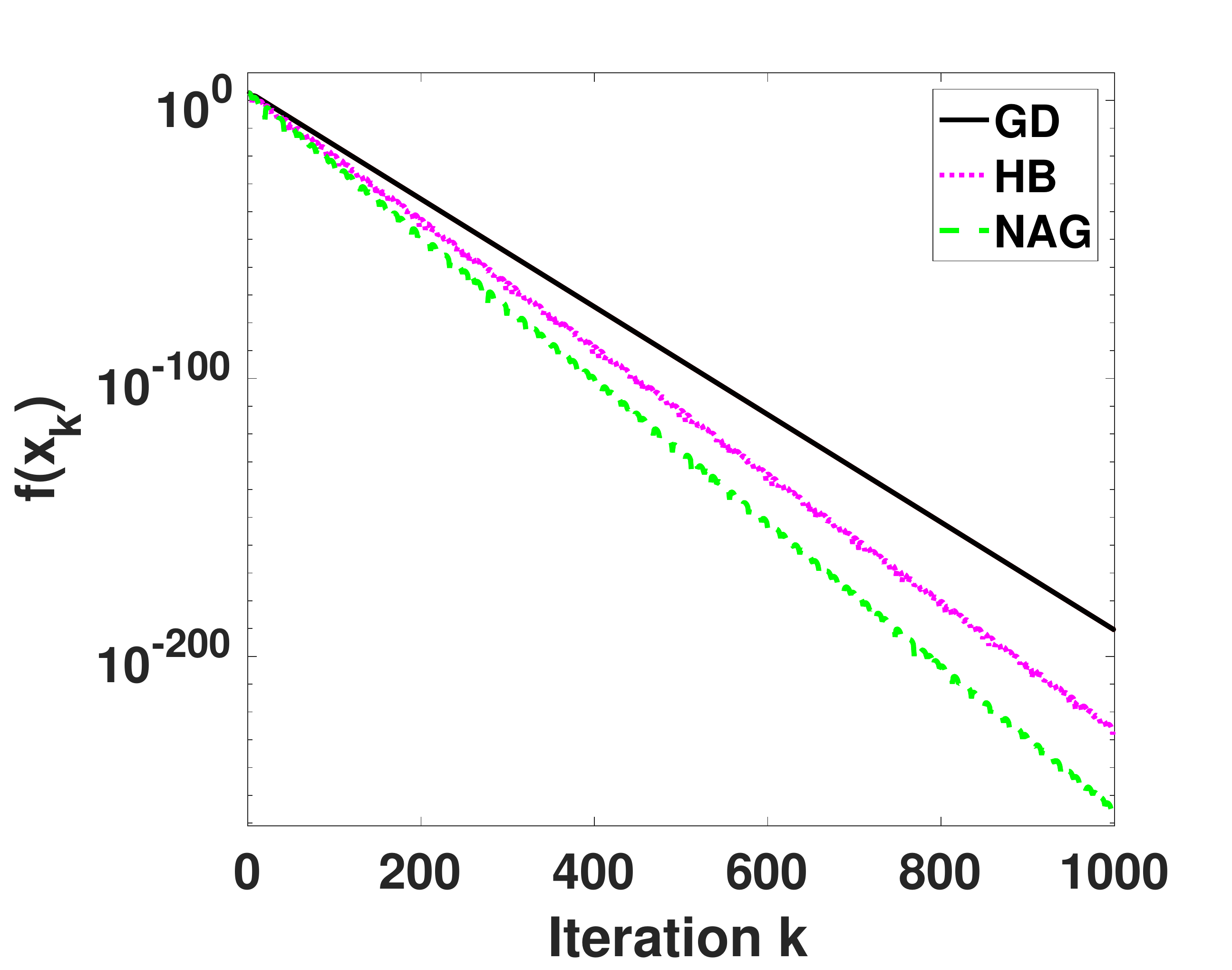} }}%
	\caption{\small An example satisfying RC and numerical experiments}%
	\label{Fig:numRes}
\end{figure*}

Consider a loss function as shown in Figure~\ref{Fig:numRes}(a) with an expression as
$$f(x)=\left\{\begin{aligned}
&x^2, \quad x\in [-6,6]\\
&x^2 + 1.5|x|\left( \cos(|x|-6) -1 \right),\quad \text{otherwise}.
\end{aligned}
\right.
$$
This nonconvex loss function was also discussed in \cite{ChenCandes15solving,chi2019nonconvex}. One can check that $f$ satisfies $RC(0.5,0.5)$.

We initialize at $x_0=x_1=24$ and choose $\alpha=0.1$. By Theorem \ref{thm:HBregionInf}, HB can converge when $\beta<0.5942$. For NAG, we choose $\beta<0.6950$ according to Theorem \ref{thm:NesregionInf}. Furthermore, it is common to choose the hyper-parameter $\beta$ as large as possible to obtain a better performance. Therefore, the corresponding $\beta$'s to HB and NAG are chosen as $0.59$ and $0.69$, respectively. In Figure \ref{Fig:numRes}(b), we see that all the three algorithms can converge and the two accelerated methods HB and NAG obviously outperform GD.

\section{Local Regularity Condition}\label{sec:LRC}

So far, all the above derivations assume RC holds globally.
In addition, the existing control framework for optimization methods \cite{lessard2016analysis,fazlyab2018analysis,hu2017unified,hu2017dissipativity,van2018fastest,cyrus2018robust,dhingra2018proximal} all require global constraints. 
In certain problems, however, RC may only hold locally around the fixed point. In this section, we explain how our framework can accommodate such cases, as long as the algorithms are initialized properly as stated in the following theorem whose proof can be found in the appendix. 

\begin{thm}\label{thm:localIni}
	Let $x^\star\in\mathbb{R}^n$ be a local minimizer of the loss function $f(\cdot)$ which satisfies RC($ \mu,\lambda,\epsilon $) with some positive constants $\mu,\lambda,\epsilon$. Assume that $P\succ 0$ is a feasible solution of the LMI \eqref{eq:maintheo}. If the first two iterates initialized properly according to $z_{-1},z_0\in\mathcal{N}_{\epsilon/\sqrt{10\mathrm{cond}(P)}}(x^\star)$, then $y_k\in \mathcal{N}_\epsilon(x^{\star}),\forall k$.
\end{thm}

Theorem~\ref{thm:localIni} ensures all the following iterates will not exceed the local neighborhood satisfying RC since $y_k\in \mathcal{N}_\epsilon(x^{\star})$ for all $k$, so that we can still transfer RC to a quadratic bound at each iteration, and thus all the previous results still hold for the convergence analysis of AGD under a general setting where RC only holds locally.

 In practice, spectral methods can be used as an initialization scheme to locate an initial estimate in a desired neighborhood of the fixed point. For example, we consider a popular inverse problem in signal processing called phase retrieval, where the goal is to recover a signal $x^\star\in\mathbb{R}^n$ from the magnitudes of its linear measurements, $y_r=|a_r^Tx^\star|^2$, $r=1,\ldots, m$, where $a_r\in\mathbb{R}^{n}$ is the $r$th sampling vector, and $m$ is the number of samples. If $a_r$'s are drawn with i.i.d. standard Gaussian entries, it is shown that the loss function $f(z)=\frac{1}{2m}\sum_{r=1}^{m}\left(y_r - |a_r^Tz|^2 \right)^2$ satisfies RC locally in $\mathcal{N}_{\epsilon}(x^{\star})$ (ignoring the sign ambiguity in identifying $x^\star$), where $\epsilon$ is a small constant (e.g. $1/10$) \cite{candes2015phase}. On the other end, the spectral method proposed in \cite{ChenCandes15solving} returns an initial estimate $\mathcal{N}_{\epsilon}(x^{\star})$ as soon as the sample complexity $m$ is above the order of $O(n/\epsilon)$. Therefore, as long as $m=O(n)$, the spectral method can successfully land an initialization in the region satisfying RC. In addition, the quality of the initialization also impacts the iteration complexity logarithmically as suggested by \eqref{eq:linear_convergence}. We refer the readers to \cite{candes2015phase,KesMonSew2010} for more details of initialization techniques.

\section{Conclusions}

In this paper, we apply control tools to analyze the convergence of AGD (including HB and NAG) under the Regularity Condition. Our main contribution lies in the {\em analytical} characterization of the convergence regions in terms of the algorithm parameters $(\alpha,\beta)$ and the RC parameters $(\mu, \lambda)$. Such convergence results do not exist in the current literature and offer useful insights in the analysis and design of AGD for a class of nonconvex optimization problems.

\section*{Acknowledgement}
The work of H. Xiong and W. Zhang is partly supported by the National Science Foundation under grant CNS-1552838. The work of Y. Chi is supported in part by ONR under grant N00014-18-1-2142, by ARO under grant W911NF-18-1-0303, and by NSF under grants CCF-1806154 and ECCS-1818571.

\appendix

\section{RC vs sector bound}

Recall the quadratic bound for $RC(\mu,\lambda)$ has the following form:
\begin{equation}\label{rcQuaBd}
\left[ \begin{array}{c}
y_k-y_*\\u_k-u_*\end{array} 
\right]^T \left[ \begin{array}{c c}
-\lambda I_n  &I_n\\
I_n  &-\mu I_n
\end{array} 
\right] 
\left[ \begin{array}{c}
y_k-y_*\\u_k-u_*\end{array} 
\right] \geq  0.
\end{equation}
The sector bound condition in \cite[Lemma 6]{lessard2016analysis} is described by the following quadratic constraint:
\begin{equation}\label{secBd}
\left[ \begin{array}{c}
y_k-y_*\\u_k-u_*\end{array} 
\right]^T \left[ \begin{array}{c c}
-2mL I_n  &(L+m)I_n\\
(L+m)I_n  &-2 I_n
\end{array} 
\right] 
\left[ \begin{array}{c}
y_k-y_*\\u_k-u_*\end{array} 
\right] \geq  0,
\end{equation}
where $m$ and $L$ are the slopes of the lines forming the sector. For simplicity, assume $m\le L$.
By comparing \eqref{rcQuaBd} and \eqref{secBd}, we can find that the quadratic bounds corresponding to RC and the sector bound are essentially the same. Specifically, given the sector bound \eqref{secBd}, we can set $\lambda=\frac{2mL}{m+L}$ and $\mu=\frac{2}{m+L}$ , which leads to  RC in \eqref{rcQuaBd}.  Similarly, given RC in \eqref{rcQuaBd}, we can immediately obtain an equivalent sector bound condition by setting $m=\frac{1-\sqrt{1-\lambda\mu}}{\mu}$ and $L=\frac{1+\sqrt{1-\lambda\mu}}{\mu}$. One special situation is when RC holds only locally, as is the case in most applications, which is handled carefully in this paper.

\section{Proof details}
\subsection{Proof of Proposition \ref{thm:convergence}}

Assume the key LMI holds with some $P\succ 0$, i.e.,
\begin{equation}\label{eq:maintheo'}
\left[\begin{array}{cc}
A^TPA-\rho^2 P &A^TPB \\
B^TPA &B^TPB
\end{array}
\right]+  \left[\begin{array}{c c}
C & D\\
\mathbf{0}_{1\times 2} &1
\end{array}
\right]^T\left[\begin{array}{c c}
-\lambda &1\\
1 &-\mu
\end{array}
\right]
\left[\begin{array}{c c}
C & D\\
\mathbf{0}_{1\times 2} &1
\end{array}
\right]\preceq0.
\end{equation}
Multiplying \eqref{eq:maintheo'} by $\left[\begin{array}{c}
\phi_k-\phi_*\\u_k-u_*\end{array}
\right]^T$ from the left and $\left[\begin{array}{c}
\phi_k-\phi_*\\u_k-u_*\end{array}
\right]$ from the right respectively, and inserting the Kronecker product, we obtain
	$$\begin{aligned}
	&\left[\! \begin{array}{c}
	\phi_k-\phi_*\\u_k-u_*\end{array} \!
	\right]^T\left(\left[\! \begin{array}{cc}
	A^TPA-\rho^2 P &A^TPB \\
	B^TPA &B^TPB
	\end{array} \!
	\right] \otimes  I_n\right)
	\left[ \begin{array}{c}
	\phi_k-\phi_*\\u_k-u_*\end{array} 
	\right]\! \\
	&+\! \left[\! \begin{array}{c}
	y_k-y_*\\u_k-u_*\end{array} \!
	\right]^T \left[\! \begin{array}{c c}
	-\lambda I_n  &I_n\\
	I_n  &-\mu I_n
	\end{array} \!
	\right] 
	\left[\! \begin{array}{c}
	y_k-y_*\\u_k-u_*\end{array} \!
	\right] \leq  0.
	\end{aligned}$$

We have that RC can be equivalently represented as a quadratic bound of the feedback term $u_k = \nabla f(y_k)$ as:
\begin{equation}
\left[\begin{array}{c}
y_k-y_*\\u_k-u_*\end{array}
\right]^T\left[\begin{array}{c c}
-\lambda I_n &I_n\\
I_n &-\mu I_n
\end{array}
\right]
\left[\begin{array}{c}
y_k-y_*\\u_k-u_*\end{array}
\right]\geq 0,
\end{equation}
which further implies that
\begin{equation}\label{LyDecay}
\left[\begin{array}{c}
\phi_k-\phi_*\\u_k-u_*\end{array}
\right]^T\left(\left[\begin{array}{cc}
A^TPA-\rho^2 P &A^TPB \\
B^TPA &B^TPB
\end{array}
\right]\otimes I_n\right) \cdot \left[\begin{array}{c}
\phi_k-\phi_*\\u_k-u_*\end{array}
\right]\leq 0.
\end{equation}

Observe $\phi_{k+1}=(A\otimes I_n) \phi_k+(B\otimes I_n) u_k$. Hence we can further rearrange and simplify \eqref{LyDecay} as
$$\begin{aligned}
&(\phi_{k+1}-\phi_*)^T \left(P\otimes I_n\right) (\phi_{k+1}-\phi_*)\\
&\leq \rho^2
(\phi_k-\phi_*)^T\left(P\otimes I_n\right)(\phi_k-\phi_*).
\end{aligned}
$$
Such exponential decay with $P\succ 0$ for all $k$ can conclude Proposition \ref{thm:convergence}.

\subsection{Proof of Lemma \ref{lem:feedbshift}}

We check the inner product of the input $u_k$ and output $y_k$ (recall that $y_*=x^{\star}$) of the dynamical system \eqref{systemshift}:
\begin{align*}
\langle u_k-u_*, y_k-y_*\rangle =& \langle -\alpha\nabla  f(y_k) - \delta(y_k-y_*), y_k-y_*\rangle\\
=& -\delta\| y_k-y_*\|^2 - \alpha\langle \nabla  f(y_k), y_k-y_*\rangle\\
\leq& -\delta\| y_k-y_*\|^2- \frac{\alpha\mu}{2}\|\nabla  f(y_k)\|^2 - \frac{\alpha\lambda}{2}\|y_k-y_*\|^2\\
=&-\frac{\mu}{2\alpha}\|u_k-u_*\|^2-\frac{\delta\mu}{\alpha}\langle u_k-u_*, y_k-y_*\rangle\\ &-\left(\delta+\frac{\alpha\lambda}{2}+\frac{\delta^2\mu}{2\alpha}\right)\|y_k-y_*\|^2.
\end{align*}
By rearrangement, we have 
$$  
-\left(2\alpha\delta+\lambda\alpha^2+\mu\delta^2\right)\|y_k-y_*\|^2
-2(\alpha+\mu\delta)\langle u_k-u_*, y_k-y_*\rangle-\mu\|u_k-u_*\|^2\geq 0,
$$
and thus conclude \eqref{RCshift}.

\subsection{Proof of Lemma \ref{lem:kypcalg}}

We check the KYPC($A',M'$) as listed for AGD \eqref{systemshift} characterized by $\alpha,\beta_1,\beta_2$:
\begin{enumerate}
	\item $\text{det}(e^{j\omega}I-A')\neq 0 $ for $\omega\in\mathbb{R}$;
	\item $A'$ is Schur stable;
	\item The left upper corner of $M'$ in \eqref{RCshift} is PSD.
\end{enumerate}

\textbf{Condition (1)}: Write
$$ \begin{aligned}
&\text{det}(e^{j\omega}I\!-\!A')\!
=\!\left\lvert\!
\begin{array}{cc}
e^{j\omega}\!-\!\left(1+\beta_1+\delta(1+\beta_2)\right) &\beta_1+\delta\beta_2\\
-1 &e^{j\omega}\\
\end{array}
\!\right\rvert\\
&=\cos^2\omega-\sin^2\omega-\left(1+\beta_1+\delta(1+\beta_2)\right)\cos\omega+\beta_1\\
&+\delta\beta_2+j(2\sin\omega\cos\omega-\left(1+\beta_1+\delta(1+\beta_2)\right)\sin\omega).
\end{aligned}
$$
By means of the opposite direction, let $\text{det}(e^{j\omega}I-A')=0$. Then we have
$$\left\{\begin{aligned}
&\cos^2\omega-\sin^2\omega-\left(1+\beta_1+\delta(1+\beta_2)\right)\cos\omega+\beta_1+\delta\beta_2=0, \\
&2\sin\omega\cos\omega-\left(1+\beta_1+\delta(1+\beta_2)\right)\sin\omega=0.
\end{aligned}
\right.
$$
From the second equality, we have $\sin\omega=0$ or $\left(1+\beta_1+\delta(1+\beta_2)\right)=2\cos\omega$, which we discuss separately. 
(a) If $\sin\omega=0$, $\cos \omega=\pm 1$. Then the first equality becomes $\left(1+\beta_1+\delta(1+\beta_2)\right)=\pm(1+\beta_1+\delta\beta_2)$; (b) if $\left(1+\beta_1+\delta(1+\beta_2)\right)=2\cos\omega$, then from the first equality we need $\cos^2\omega-\sin^2\omega-\left(1+\beta_1+\delta(1+\beta_2)\right)\cos\omega+\beta_1+\delta\beta_2=-1+\beta_1+\delta\beta_2=0$.

To conclude, condition (1) is satisfied if and only if:
\begin{equation}\label{cond2}
\left(1\!+\!\beta_1\!+\!\delta(1\!+\!\beta_2)\right)\!\neq\!\pm(1\!+\!\beta_1\!+\!\delta\beta_2),\! \beta_1\!+\!\delta\beta_2\neq 1.
\end{equation}

\textbf{Condition (2)}: We want 
$$ A'=\left[\begin{array}{cc}
\left(1+\beta_1+\delta(1+\beta_2)\right) &-(\beta_1+\delta\beta_2)\\
1 &0
\end{array}
\right]$$ 
to be Schur stable, for which it suffices to check its eigenvalues are bounded by $1$ in magnitude. We start by writing out the characteristic equation of $A'$:
$$\begin{aligned}
\lvert \lambda I-A' \rvert &=\left\lvert
\begin{array}{cc}
\lambda-\left(1+\beta_1+\delta(1+\beta_2)\right) &\beta_1+\delta\beta_2\\
-1 &\lambda\\
\end{array}
\right\rvert\\
& = \lambda^2-\lambda \left(1+\beta_1+\delta(1+\beta_2)\right)+\beta_1+\delta\beta_2=0.
\end{aligned}$$
The eigenvalues of $A'$ are given as the two roots of the above polynomial:
\begin{align*}
\lambda_{1,2} &=\frac{\left(1+\beta_1+\delta(1+\beta_2)\right)}{2} \pm \frac{\sqrt{\left(1+\beta_1+\delta(1+\beta_2)\right)^{2}-4(\beta_1+\delta\beta_2)}}{2}.
\end{align*}
We need to make sure $\lvert \lambda_{1,2} \rvert<1$. When the eigenvalues are complex-valued, i.e. $\left(1+\beta_1+\delta(1+\beta_2)\right)^{2}<4(\beta_1+\delta\beta_2)$, $\lvert \lambda_1 \rvert=\lvert \lambda_2 \rvert=\beta_1+\delta\beta_2$, and thus we need $\beta_1+\delta\beta_2<1$ in this case. When the eigenvalues are real-valued, i.e. $\left(1+\beta_1+\delta(1+\beta_2)\right)^{2}\geq 4(\beta_1+\delta\beta_2)$, 
we have
$$
4(\beta_1+\delta\beta_2)\leq \left(1+\beta_1+\delta(1+\beta_2)\right)^2<(1+\beta_1+\delta\beta_2)^2.
$$
To conclude, condition (2) is satisfied if and only if
\begin{equation}\label{cond3}
-\frac{2(1+\beta_1)}{1+2\beta_2}<\delta<0.
\end{equation}

\textbf{Condition (3)}: We want the left upper corner of $M'$ in \eqref{RCshift} to be PSD, which means 
$$-\left(2\alpha\delta+\lambda\alpha^2+\mu\delta^2\right)\geq0.$$ 
Then condition (3) is satisfied if and only if
\begin{equation}\label{cond4}
\frac{-\delta-|\delta|\sqrt{1-\mu\lambda}}{\lambda}\leq\alpha\leq\frac{-\delta+ |\delta|\sqrt{1-\mu\lambda}}{\lambda}.
\end{equation}
To conclude, by unifying all conditions \eqref{cond2}\eqref{cond3}\eqref{cond4}, we conclude that the parameters of AGD need to satisfy 
\begin{equation}
0<\alpha<\frac{2(1+\beta_1)(1+\sqrt{1-\mu\lambda})}{\lambda(1+2\beta_2)}.
\end{equation}

\subsection{Proof of Theorem \ref{thm:HBregionInf}}

For HB, we set $\beta_1=\beta$, $\beta_2=0$ and denote $u:=\cos\omega$, then \eqref{FDIgeneral} can be rewritten as:
\begin{equation}\label{FDI}
 h(u)\!:=\!4\mu\beta u^2\!-\!2(2\mu\beta\!+\!\mu\beta^2\!-\!\alpha\beta\!+\!\mu\!-\!\alpha)u + 2\mu\!+\!2\mu\beta^2\!-\!2\alpha\!-\!2\alpha\beta\!+\!\lambda\alpha^2\!\geq\! 0, \quad \forall u\in[-1,1]. 
\end{equation}
Observe that $h(u)$ is a quadratic function depending on $u$. We can check the minimal value of $h(\cdot)$ on $[-1,1]$ by discussing its axis of symmetry denoted as $S=\frac{2\mu\beta+\mu\beta^2-\alpha\beta+\mu-\alpha}{4\mu\beta}$.
\begin{enumerate}
	\item When $S\geq 1$, $\alpha\leq\frac{\mu(1-\beta)^2}{1+\beta}$. Then
	$$h(u)_{\min}=h(1)=\lambda\alpha^2>0.
	$$
	Thus the feasible region in this case is:
	\begin{equation}\label{region1}
	\left\{(\alpha,\beta):\alpha\leq\frac{\mu (\beta-1)^2}{\beta+1}, 0<\beta<1 \right\}.
	\end{equation}
	\item When $S\leq -1$, $\alpha\geq\frac{\mu\beta^2+6\mu\beta+\mu}{1+\beta}$. We need
	$$h(u)_{\min}=h(-1)=\lambda\alpha^2-4(1+\beta)\alpha+4\mu(1+\beta)^2\geq 0.
	$$
	Thus the feasible region in this case is:
	\begin{equation}\label{region2}
	\hspace{-1.3em}\begin{aligned}
	&\left\{\!(\alpha\!,\beta)\!:\!\alpha\geq\frac{2(\beta\!+\!1)(1\!+\!\sqrt{1\!-\!\mu\lambda})}{\lambda}, 0\!<\!\beta\!<\!1\! \right\}\!\\ 
	&\begin{aligned}
	\cup \left\{\!(\alpha\!,\beta)\!:\!\frac{\mu \beta^2+6\mu\beta+\mu}{\beta+1}\leq\alpha \leq\frac{2(\beta+1)(1-\sqrt{1-\mu\lambda})}{\lambda}, 0<\beta<1 \right\}.
	\end{aligned}
	\end{aligned}
	\end{equation}
	\item When $-1<S<1$, $\frac{\mu\beta^2+6\mu\beta+\mu}{1+\beta}<\alpha<\frac{\mu(1-\beta)^2}{1+\beta}$. We want $
	h(u)_{\min} =h\left(\frac{2\mu\beta+\mu\beta^2-\alpha\beta+\mu-\alpha}{4\mu\beta}\right)\geq 0$. It is equivalent to solve 
	$$\hspace{-0.5em}
	\begin{aligned}
	&(4\mu\lambda\beta-\beta^2-1-2\beta)\alpha^2-\left(2\mu\beta+2\mu\beta^2-2\mu\beta^3\right.\\
	&\left.-2\mu\right)\alpha+4\mu^2\beta^3+4\mu^2\beta-6\mu^2\beta^2-\mu^2\beta^4-\mu^2\geq 0.
	\end{aligned} $$
	Since $4\mu\lambda\beta-\beta^2-1-2\beta<4\beta-\beta^2-1-2\beta\leq0$,
	Thus the feasible region in this case is:
	\begin{equation}\label{region3}
	 \left\{\! (\alpha,\beta):\!\frac{\mu (\beta-1)^2}{\beta+1}\!\leq\!\alpha\!\leq\!R, 0\!<\!\beta\!<\!1\!\right\}  \cap 
	\left\{\! (\alpha,\beta):\!\alpha\!\leq\!\frac{\mu \beta^2\!+\!6\mu\beta\!+\!\mu}{\beta+1} \!, 0\!<\!\beta\!<\!1\! \right\} ,
	\end{equation}
	where $R=\frac{P_2 - \sqrt{P_2^2-4P_1P_3}}{2P_1},
	P_1=4\mu\lambda\beta-\beta^2-1-2\beta, P_2= 2\mu\beta+2\mu\beta^2-2\mu\beta^3-2\mu, P_3=4\mu^2\beta^3+4\mu^2\beta-6\mu^2\beta^2-\mu^2\beta^4-\mu^2 $.
\end{enumerate}

Taking the union of \eqref{region1}\eqref{region2}\eqref{region3} gives the result of the FDI \eqref{FDI}. Further intersecting with the condition \eqref{HbshiftLem} obtained in Lemma \ref{lem:kypcalg} leads to the final region.

\subsection{Proof of Theorem \ref{thm:NesregionInf}}

Similar with the proof of Theorem \ref{thm:HBregionInf}, for NAG, we set $\beta_1=\beta_2=\beta$ and denote $u:=\cos\omega$, then \eqref{FDIgeneral} can be rewritten as:
\begin{equation}\label{FDI_Nes}
\hspace{-0.3em}\begin{aligned}
h(u)& := (-4\mu\beta+4\alpha\beta)u^2 +\left[2(\mu-\alpha)(1+\beta)^2+2(\lambda\alpha-1)\alpha(1+\beta)\beta \right]u\\
&-2\mu(1+\beta^2)+2\alpha(1+\beta)^2-\lambda\alpha^2\beta^2 -2\alpha(1-\beta)\beta-\lambda\alpha^2(1+\beta)^2\leq 0, \quad  \forall u\in[-1,1].
\end{aligned}	
\end{equation}
We check the maximal value of the quadratic function $h(\cdot)$ on $[-1,1]$.

\begin{enumerate}
	\item When $\alpha=\mu$, $h(u)_{max}=f(-1)\leq 0$. Thus the feasible region in this case is:
	\begin{equation}\label{regionNes1}
	\hspace{-0.55em}\left\{\!(\alpha,\beta)\!:\alpha\!=\!\mu, 0\!<\!\beta\!\leq\!\frac{-1\!+\!\mu\lambda\!+\!\sqrt{1\!-\!\mu\lambda}}{2(1-\mu\lambda)} \right\}.
	\end{equation}
	\item When $\alpha>\mu$, we need to let $h(1)\leq 0,h(-1)\leq 0$. Since $h(1)=-\lambda\alpha^2<0$, we only need to check $h(-1)$. Thus the feasible region in this case is:
	\begin{equation}\label{regionNes2}
	\left\{(\alpha,\beta):\alpha\geq L_1, 0<\beta< 1 \right\} \cup
	\left\{(\alpha,\beta):\mu<\alpha\leq R_1, 0<\beta<1 \right\},
	\end{equation}
	where $L_1=\frac{2(\beta+1)(1+\sqrt{1-\mu\lambda})}{\lambda(1+2\beta)},R_1=\frac{2(\beta\!+\!1)(1\!-\!\sqrt{1\!-\!\mu\lambda})}{\lambda(1+2\beta)}$.
	\item When $\alpha<\mu$, we can check the maximal value of the quadratic function $h(u)$ by discussing the axis of symmetry $S = \frac{(\mu-\alpha)(1+\beta)^2+(\lambda\alpha-1)\alpha(1+\beta)\beta}{4\mu\beta-4\alpha\beta}$.
	\begin{enumerate}
		\item When $S\geq 1$, $h(u)_{max}=-\lambda\alpha^2<0$. Thus the feasible region in this case is:
		\begin{equation}\label{regionNes3}
		\left\{(\alpha,\beta):\alpha\!\leq\!\frac{B_1\!-\!\sqrt{B_1^2\!-\!C_1}}{2\lambda\beta(\beta+1)}, 0\!<\!\beta\!<\!1\! \right\},
		\end{equation}
		where 
		$B_1=1-\beta+2\beta^2,C_1=4\mu\lambda\beta(1+\beta)(1-\beta)^2.$
		\item When $S\leq -1$, $h(u)_{max}=h(-1)\leq 0$. Thus the feasible region in this case is:
		\begin{equation}\label{regionNes4}
		\hspace*{-0.3cm}
		\left\{(\alpha,\beta):\frac{B_2-\!\sqrt{B_2^2\!-\!C_2}}{2\lambda\beta(\beta+1)}\!\leq\!\alpha\!<\!\mu, 0\!<\!\beta\!<\!1\! \right\},
		\end{equation}
		where
		$B_2=1\!+\!7\beta\!+\!2\beta^2,C_2=4\mu\lambda\beta(1+\beta)(1+6\beta+\beta^2).$
		\item When $-1<S<1$, $h(u)_{max}=h(S)\leq 0$. Thus the feasible region in this case is:
		\begin{equation}\label{regionNes5}
		\hspace{-0.58em}\left\{(\alpha,\beta)\!:L\!<\!\alpha\!<\!R, 0\!<\!\beta\!<\!1, g(S)\!\leq\! 0 \right\},\!
		\end{equation}
		where $L= \frac{1-\beta+2\beta^2\!-\!\sqrt{(1-\beta+2\beta^2)^2-4\mu\lambda\beta(1+\beta)(1-\beta)^2}}{2\lambda\beta(\beta+1)}$,
		$R = \frac{1+7\beta+2\beta^2\!-\!\sqrt{(1+7\beta+2\beta^2)^2-4\mu\lambda\beta(1+\beta)(1+6\beta+\beta^2)}}{2\lambda\beta(\beta+1)}$.
		and $g(S)$ is noted in Remark \ref{rem:Nes}, that is,
		$$\begin{aligned}
		g(S):=&4\mu\beta S^2-2(2\mu\beta+\mu\beta^2-\alpha\beta+\mu-\alpha)S\\
		&+2\mu+2\mu\beta^2-2\alpha-2\alpha\beta+\lambda\alpha^2.
		\end{aligned}
		$$
	\end{enumerate}
\end{enumerate}

The result of the FDI \eqref{FDI_Nes} is the union of all the above regions \eqref{regionNes1}\eqref{regionNes2}\eqref{regionNes3}\eqref{regionNes4}\eqref{regionNes5}. 
Further intersecting with the condition \eqref{HbshiftLem} obtained in Lemma \ref{lem:kypcalg} leads to the final region.

\subsection{Proof of Theorem \ref{thm:localIni}}

Since $z_{-1},z_0\in\mathcal{N}_{\epsilon/\sqrt{10\mathrm{cond}(P)}}({x^\star})$, we have $\|\phi_0-\phi_*\|<\frac{\epsilon}{\sqrt{5\mathrm{cond}(P)}}$. The exponential decay with $P\succ 0$: $(\phi_{k+1}-\phi_*)^TP(\phi_{k+1}-\phi_*)\leq \rho^2
(\phi_k-\phi_*)^TP(\phi_k-\phi_*)
$ implies that 
$ \|\phi_{k}-\phi_*\|\leq \sqrt{\mathrm{cond}(P)}\rho^k\|\phi_{0}-\phi_*\|
$, which we have argued in Subsection \ref{subsec:lyFunc}. Therefore,
$$\begin{aligned}
\|\phi_k-\phi_*\|&\leq \sqrt{\mathrm{cond}(P)}\rho^k\|\phi_{0}-\phi_*\|\\
&<\sqrt{\mathrm{cond}(P)} \cdot \|\phi_{0}-\phi_*\|\\
&<\sqrt{\mathrm{cond}(P)}\cdot\frac{\epsilon}{\sqrt{5\mathrm{cond}(P)}}\\
&<\epsilon/\sqrt{5}.
\end{aligned}
$$
As a consequence, 
$$\begin{aligned}
\|y_k-x^\star\|&=\left\|(1+\beta_2)z_{k}^{(1)}-\beta_2 z_{k}^{(2)}-x^\star\right\|\\
&=\left\|C (\phi_k-\phi_*)\right\|\\
&\leq \left\|C^T\right\| \left\|\phi_k-\phi_*\right\|\\
&< \sqrt{(1+\beta_2)^2+\beta_2^2}\cdot\left(\epsilon/\sqrt{5}\right) \\
&<\epsilon,
\end{aligned}
$$
where we recall $C=[ 1+\beta_2 ,-\beta_2 ]$, and the last line used $\beta_2<1$.

\bibliographystyle{IEEEtran}
\bibliography{bibfileNonconvex2018}

\end{document}